\documentclass[a4paper]{amsart}
\usepackage{fixltx2e}
\usepackage[utf8]{inputenc}
\usepackage[T1]{fontenc}
\usepackage{amssymb,amsthm,amsmath}
\usepackage[british]{babel}
\usepackage[all]{xy}
\usepackage{calrsfs}
\usepackage{hyperref}
\usepackage{mathbbol}
\usepackage{mathabx}

\theoremstyle{plain}
\newtheorem*{theorem}{Theorem}

\newcommand{\defn}{\textbf}

\DeclareMathOperator{\Kernel}{Ker}
\newcommand{\Ker}{\Kernel}
\newcommand{\C}{\ensuremath{\mathcal{C}}}
\newcommand{\N}{\ensuremath{\mathbb{N}}}
\newcommand{\Mon}{\ensuremath{\mathsf{Mon}}}

\newdir{>>}{{}*!/3.5pt/:(1,-.2)@^{>}*!/3.5pt/:(1,+.2)@_{>}*!/7pt/:(1,-.2)@^{>}*!/7pt/:(1,+.2)@_{>}}
\newdir{ >>}{{}*!/8pt/@{|}*!/3.5pt/:(1,-.2)@^{>}*!/3.5pt/:(1,+.2)@_{>}}
\newdir{ |>}{{}*!/-3.5pt/@{|}*!/-8pt/:(1,-.2)@^{>}*!/-8pt/:(1,+.2)@_{>}}
\newdir{ >}{{}*!/-8pt/@{>}}
\newdir{>}{{}*:(1,-.2)@^{>}*:(1,+.2)@_{>}}
\newdir{<}{{}*:(1,+.2)@^{<}*:(1,-.2)@_{<}}

\begin{document}

\title[A new characterisation of groups amongst monoids]{A new characterisation of\\ groups amongst monoids}
\author{Xabier García-Martínez}
\address[Xabier García-Martínez]{Department of Algebra, University of Santiago de Compostela, 15782 Santiago de Compostela, Spain.}
\email{xabier.garcia@usc.es}

\thanks{
	Supported by Ministerio de Economía y 
	Competitividad (Spain), grant MTM2013-43687-P (European 
	FEDER support included), by Xunta de Galicia, 
	grant GRC2013-045 (European FEDER support included), and 
	by an FPU scholarship of the Ministerio de Educación, Cultura y Deporte (Spain). The author would like to thank Dr.~Tim Van der Linden for his help and dedication, and the Institut de Recherche en Mathématique et Physique (IRMP) for its kind hospitality during his stays in Louvain-la-Neuve}

\begin{abstract}
We prove that a monoid $M$ is a group if and only if, in the category of monoids, all points over $M$ are strong. This sharpens and greatly simplifies a result of Montoli, Rodelo and Van der Linden~\cite{MRVdL-TCOGAM} which characterises groups amongst monoids as the protomodular objects.
\end{abstract}


\maketitle

In their article~\cite{MRVdL-TCOGAM}, Montoli, Rodelo and Van der Linden introduce, amongst other things, the concept of a \emph{protomodular object} in a finitely complete category~$\C$ as an object $Y\in \C$ over which all points are \emph{stably strong}. The aim of their definition is two-fold: first of all, to provide a categorical-algebraic characterisation of groups amongst monoids as the protomodular objects in the category $\Mon$ of monoids; and secondly, to establish an object-wise approach to certain important conditions occurring in categorial algebra such as protomodularity~\cite{Bourn1991, Borceux-Bourn} and the Mal'tsev axiom~\cite{CLP, CPP}. 

We briefly recall some basic definitions; see~\cite{Bourn-monad, MartinsMontoliSobral2, MRVdL-TCOGAM} for more details. Let $\C$ be a finitely complete category, which we also take to be pointed for the sake of simplicity. In $\C$, a pair of arrows $(r\colon {W\to X}, s\colon {Y\to X})$ is  \defn{jointly strongly epimorphic} when if $mr'=r$, $ms'=s$ for some given monomorphism $m\colon {M\to X}$ and arrows $r'\colon{W\to M}$, $s'\colon{Y\to M}$, then $m$ is an isomorphism. In the case of monoids, this means that any $x\in X$ can be written as a product $r(w_{1})s(y_{1})\cdots r(w_{n})s(y_{n})$ for some $w_{j}\in W$, $y_{j}\in Y$. This characterisation follows easily from the fact that $(r,s)$ is a jointly strongly epimorphic pair in $\Mon$ if and only if the induced monoid morphism ${W+Y\to X}$ is a surjection---see, for instance, \cite[Corollary~A.5.4 combined with Example~A.5.16]{Borceux-Bourn}. Given an object $Y$ in $\C$, a \defn{point over $Y$} is a pair of morphisms $(f\colon{X\to Y},s\colon{Y\to X})$ such that $fs=1_{Y}$. A~point $(f,s)$ is said to be \defn{strong} when the pair $(\ker(f)\colon{\Ker(f)\to X}, s\colon {Y\to X})$ is jointly strongly epimorphic. The point $(f,s)$ is \defn{stably strong} when all of its pullbacks are strong. More precisely, if $g\colon {Z\to Y}$ is any morphism, then the pullback $g^{*}(f)$ together with its splitting induced by $s$ is a strong point.

Even though the concept of a protomodular object serves the intended purpose of characterising groups amongst monoids, the proof of this characterisation given in~\cite{MRVdL-TCOGAM} is rather complicated, since it relies on another, more subtle, characterisation in terms of the so-called \emph{Mal'tsev objects}. The present short note aims to improve the situation by giving a quick and direct proof of a more general result: a monoid is a group as soon as all points over it are strong.

\begin{theorem}
A monoid $M$ is a group if and only if, in $\Mon$, all points over $M$ are strong.
\end{theorem}
\begin{proof}
It is shown in~\cite{SchreierBook}---this is Proposition~2.2.4 combined with Lemma~2.1.6---that for any group $M$, all points over it are \emph{homogenous}, which makes them (stably) strong. So in particular, if $M$ is a group, then all points over $M$ are strong. We prove the other implication.

Consider $m\in M$ and the induced split extension
\[
\xymatrix@C=4em{ 0 \ar[r] & K \ar@{{ |>}->}[r] & \N + M \ar@{-{>>}}@<-.5ex>[r]_-{\lgroup m\; 1_{M}\rgroup} & M \ar@{{ >}->}@<-.5ex>[l]_-{\iota_{M}} \ar[r] & 0,}
\]
where $m\colon {\N\to M}$ is the morphism which sends the generator $1$ of $(\N,+,0)$ to the element $m$ of $M$. By the assumption that $(\lgroup m\; 1_{M}\rgroup,\iota_{M})$ is a strong point, $1\in \N$ can be written as
\[
1=k_{1}m_{1}\cdots m_{i}k_{i+1}m_{i+1}\cdots k_{n}m_{n}
\]
for some $k_{j}\in K$ and $m_{j}\in M$. Since $1$ is not invertible in $\N$, it must appear in exactly one of the factors $k_{j}$ in the product on the right, say in $k_{i+1}$. Then neither $k_{1}m_{1}\cdots m_{i}$ nor $m_{i+1}\cdots k_{n}m_{n}$ contains any non-zero elements of $\N$, so we have that in $\N + M$
\[
1 = a'kb'
\]
for some $a'$, $b'\in M$ and $k\in K$. Since $1$ appears in $k$ we can write $k=a1b$ where $a$, $b\in M$. Necessarily then $e_{M}=a'a$ and $e_{M}=bb'$, because $1=a'a\cdot 1\cdot bb'$. Furthermore, since $k$ is in the kernel of $\lgroup m\; 1_{M}\rgroup$, we also have that $e_{M}=amb$. So, clearly, $a$ and $b$ are invertible. As a consequence, $m$ is invertible as well. We conclude that $M$ is a group.
\end{proof}

Note that the above proof shows in particular why $M $ is \defn{gregarious} in the sense of~\cite{Borceux-Bourn}, which means that for any $m$ there exist $a$ and~$b$ such that $e_{M}=amb$. However, the proof also shows that those $a$ and~$b$ are invertible, and thus $M$ is a group.

This result seems to indicate that in certain cases (like, for instance, in the category of monoids) it makes sense to weaken the definition of a protomodular object $M$---all points over $M$ are stably strong---to the condition that those points are strong. This, and related considerations, will be the subject of future joint work with the authors of~\cite{MRVdL-TCOGAM}.


\providecommand{\noopsort}[1]{}
\providecommand{\bysame}{\leavevmode\hbox to3em{\hrulefill}\thinspace}
\providecommand{\MR}{\relax\ifhmode\unskip\space\fi MR }
\providecommand{\MRhref}[2]{%
  \href{http://www.ams.org/mathscinet-getitem?mr=#1}{#2}
}
\providecommand{\href}[2]{#2}

\end{document}